\documentclass[10pt,conference]{IEEEtran}
\usepackage{graphicx}

\usepackage{cite}
\ifCLASSINFOpdf
\else
\fi

\usepackage{mathrsfs,amsfonts,amssymb,amsthm}
\usepackage[cmex10,intlimits]{amsmath}
\usepackage{array}

\usepackage{mdwmath}
\usepackage{mdwtab}

\begin{document}
\newtheorem{theorem}{Theorem}
\newtheorem{proposition}{Proposition}
\newtheorem{corollary}{Corollary}
\newtheorem{lemma}{Lemma}
\newtheorem{definition}{Definition}
\newtheorem{assumption}{Assumption}
\newtheorem{remark}{Remark}
\newtheorem{operation}{Operation}
\newcommand{\ams}{\textit{Ann. Math. Statist.}}
\newcommand{\dfn}{\stackrel{\triangle}{=}}
\newcommand{\argmax}{\operatornamewithlimits{arg\,max}}
\newcommand{\argmin}{\operatornamewithlimits{arg\,min}}
\newcommand{\argsup}{\operatornamewithlimits{arg\,sup}}
\newcommand{\arginf}{\operatornamewithlimits{arg\,inf}}

\title{Decentralized Multihypothesis Sequential Detection}

\author{\IEEEauthorblockN{Yan Wang and Yajun Mei}
\authorblockA{School of Industrial and Systems Engineering\\
Georgia Institute of Technology,
Atlanta, Georgia 30332--0225\\
Email: \{ywang67, ymei\}@isye.gatech.edu}
}

\maketitle

\begin{abstract}
\boldmath This article is concerned with decentralized sequential testing of multiple hypotheses. In a sensor network system with limited local memory, raw observations are observed at the local sensors, and quantized into binary sensor messages that are sent to a fusion center, which makes a final decision. It is assumed that the raw sensor observations are distributed according to a set of $M \ge 2$ specified distributions, and the fusion center has to utilize quantized sensor messages to decide which one is the true distribution. Asymptotically Bayes tests are offered for decentralized multihypothesis sequential detection by combining three existing methodologies together: tandem quantizers, unambiguous likelihood quantizers, and randomized quantizers.
\end{abstract}

\IEEEpeerreviewmaketitle

\section{Introduction}\label{sec:intro}

As a subfield of signal detection or hypothesis testing, multihypothesis sequential detection has many important engineering applications such as target detection in multiple-resolution radar, serial acquisition of direct-sequence spread spectrum signals and fault detection, see Baum and Veeravalli \cite{bv}. The centralized version has been studied in both statistical and engineering literature, see the award winning papers by Dragalin, Tartakovsky and Veeravalli \cite{dtv}, \cite{dtv2} and their references for the latest development.

In recent decades the decentralized version of signal detection or hypothesis testing  has gained a great deal of attention, partly because geographically distributed sensors have been employed into a wide range of areas like military surveillance \cite{ts}, target tracking and classification \cite{lwhs}, and data filtering \cite{yllz}, etc. In the decentralized version, it is standard to assume that raw observations are observed at the local sensors, and quantized into sensor messages that are sent to a fusion center, which makes a final decision.   Unfortunately, most research on decentralized detection deals with the off-line setting and  research for the online or {\it sequential} setting is rather limited. To the best of our knowledge, so far existing research on decentralized sequential detection  is restricted to two-hypothesis, see Veeravalli, Basar, and Poor \cite{vbp}, Veeravalli \cite{vee} and Mei \cite{mei}.

The goal of this paper is to develop asymptotic optimality theory for decentralized sequential detection when there are $M \ge 2$ possible hypotheses on the models of the sensor network system. A main challenge is how to find good quantizers at the local sensors so that the fusion center is able to utilize quantized sensor messages to make effective decisions. Intuitively, the choice of good quantizers should depend on the true unknown distribution of raw sensor observations. Since there are $M \ge 2$ hypotheses, it is expected that stationary quantizers will not lead to (asymptotically) optimal tests no matter how clever one chooses it. It turns out that by combining three existing methodologies together: ``tandem quantizers" in  Mei \cite{mei}, ``unambiguous likelihood quantizers" (ULQ) in Tsitsiklis \cite{tsi}, and randomized quantizers, we are able to find good quantizers and use them to offer a family of asymptotically Bayes tests for decentralized multihypothesis sequential detection.

The remainder of this article is organized as follows. Section II provides a formal mathematical formulation of decentralized sequential multihypothesis testing problem and introduce the notation of randomized quantizer. Section III discusses tandem quantizers and constructs a family of ``two-stage'' decentralized sequential tests. This leads to a natural definition of ``maximin quantizers," in which the corresponding two-stage decentralized sequential tests are shown to be asymptotically Bayes. In Section IV, the maximin quantizers are characterized in more details as a randomized  quantizer based on at most $M-1$ (deterministic) ULQs, and numerical algorithms are provided to solve them explicitly. Section V provides specific examples to illustrate the method developed in previous sections.

\section{Notations and Problem Formulation}

\begin{figure}
\centering
\includegraphics[width=3in]{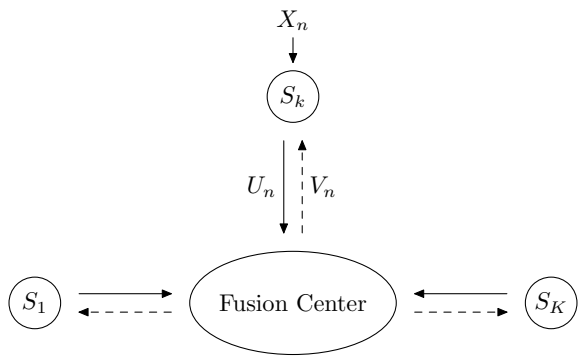}
\caption{A Sensor Network}
\label{fig:snsrnet}
\end{figure}

Fig. \ref{fig:snsrnet} shows a widely used configuration of sensor networks, where a fusion center is associated with a set of remote local sensors $S_1,\dots,S_K.$ To highlight our main ideas, we assume $K=1$ here, since the extension to systems with multiple  sensors is relatively straightforward as long as the sensor observations are independent from sensor to sensor conditioned on each hypothesis. The local sensor  takes a sequence of independent and identically distributed (i.i.d.) raw observations $X_{1}, X_{2}, \cdots$ over time $n$.
In the decentralized version, it is assumed that  the fusion center has no direct access to the raw sensor data $X_n$'s due to communication constraints. Rather, the local sensor compresses $X_{n}$ into quantized message $U_{n}\in\{0,1,\dots, l-1\},$ and sends it to the fusion center, which will then use the $U_{n}$'s as inputs to make a final decision. For our purpose, we also assume that the fusion center can send feedbacks $V_{n-1}$ to local sensor so that the local sensor can adaptively adjust sensor policies to the optimal one. For simplicity, we further assume the quantized messages to be binary, i.e., $U_n\in\{0,1\}$.

Mathematically, at time $n,$ the sensor message $U_{n}$ and fusion center feedback $V_{n-1}$ can be defined as
\begin{equation*}\label{equ:lmtdlclMmry}
U_n=\phi_n(X_n; V_{n-1})\in\{0,1\},\quad V_{n-1}=\psi_{n}(U_{[1,n-1]}),
\end{equation*}
where $U_{[1,n-1]}=(U_1,\dots,U_{n-1})$. Note that the feedback $V_{n-1}$ should only depend on the past sensor messages. Here no restrictions are imposed on $V_{n-1}$, but it turns out that $\log_{2}(M)$-bit feedbacks will be sufficient to construct asymptotically optimal tests under our setting.

In decentralized multihypothesis sequential detection, it is assumed that there are $M \ge 2$ hypotheses regarding the true  probability distribution ${\bf P}$ of $X_n$'s:
\begin{equation}\label{equ:Hypothesis}
  \textbf{H}_m: \textbf{P} = \textbf{P}_m,
\end{equation}
for $m=0,1,\dots,M-1,$ where the $X_n$'s have a probability density (or mass) function $f_m(x)$ under $\textbf{P}_m.$
Furthermore, the sensor network system will continue taking observations until the fusion center believes that there is sufficient evidence from the quantized messages $U_{n}$'s to make a final decision. That is, at a stopping time $N,$ the fusion center makes a decision $D \in \{0,1,\dots,M-1\},$ where $\{D = m\}$ means that one accepts the hypothesis $\textbf{H}_{m}.$ Here we emphasize that  the decision $\{N=n\}$ only depends on the first $n$ sensor messages, i.e., $N$  is a stopping time with respect to the filtration $\{\mathcal{F}_n=\sigma\{U_{[1,n]}\}\}$ and $D$ is measurable to $\mathcal{F}_N$.

In summary, a decentralized sequential test $\delta$ includes a sequence of quantizers $\phi_n$ at the local sensor, a sequence of feedback functions $\psi_n$, a stopping time $N$  at the fusion center and the final decision $D$.

As in Wald \cite{wald47} and Veeravalli et al. \cite{vbp}, we consider the Bayes formulation of decentralized multihypothesis sequential detection. Let $c>0$ be the cost of data sampling per time step, and $W(m,m')$ be the loss of making decision $D=m'$ when the true state of nature is $\textbf{P}_m$. We assume that all $W(m,m')$'s are non-negative and $W(m,m')=0$ if and only if $m=m'$. Let the total risk of a test $\delta$ when the true state is $m$ be
\begin{displaymath}
  \mathcal{R}_c(\delta; m)=c\textbf{E}_mN+\sum_{m'}W(m,m')\textbf{P}_m[D=m'].
\end{displaymath}
Assigning prior probabilities $\pi=(\pi_{0},\dots,\pi_{M-1})$ to $\textbf{H}_0,\cdots,\textbf{H}_{M-1},$ define the average risk of a decentralized sequential test $\delta$ as
\begin{equation}\label{equ:ttlCstsApriorPi}
  \mathcal{R}_c(\delta)=\sum_m \pi_m \mathcal{R}_c(\delta; m).
\end{equation}
The Bayesian formulation of decentralized sequential detection problems can be stated as follows:

\medskip
\textit{Problem (P1):} Minimize the $\mathcal{R}_c(\delta)$ in (\ref{equ:ttlCstsApriorPi}) among all possible decentralized sequential hypothesis testing procedures $\delta$.
\medskip

Let $\delta_{B}^*(c)$ denote a Bayes solution to (\textit{P1}), i.e., $\delta_{B}^*(c)=\argmin_{\delta}\{\mathcal{R}_c(\delta)\}$. Unfortunately, the exact form of $\delta_B^*(c)$ is too complicated to be tractable for multihypothesis sequential detection even for the centralized version, see, for example, Dragalin, Tartakovsky and Veeravalli \cite{dtv}.
This leads us to consider the ``asymptotic optimality'' approach as follows: 

\medskip
\textit{Problem (P2)}: Find a family of decentralized sequential tests $\{\delta(c)\}$ such that
\begin{equation*}\label{equ:AsmptOptmlTest}
\lim_{c\to 0}\mathcal{R}_c(\delta_B^*(c))/\mathcal{R}_c(\delta(c))=1,
\end{equation*}
where $c$ is the unit cost in (\ref{equ:ttlCstsApriorPi}). 

Problem (P2) is meaningful in application because it is often the case that the cost of doing a round of sampling is much smaller than that of making an incorrect decision.

In the remainder of this section, let us discuss the concepts of randomized quantizers and Kullback-Leibler (K-L) divergences. Denote by $\Phi$ the set of deterministic quantizers that consists of all measurable functions from $\mathbb{R}$ to $\{0,1\}$. For a quantizer $\phi\in\Phi$, let $f_m(\cdot; \phi)$ denote the induced distribution of the quantized data $\phi(X_n)$ under $\textbf{P}_m$, i.e., for $u\in\{0,1\}$, $f_m(u; \phi)=\textbf{P}_m(\phi(X_n)=u)$. Recall that the K-L divergence of $\phi$ of any state $m$ against any other state $m'\neq m$ is defined as 
\[I(m,m';{\phi})=\sum_{u=0}^{1}f_m(u; \phi)\log \frac{f_m(u; \phi)}{f_{m'}(u; \phi)}.\]
Now define a ``randomized quantizer'' $\bar{\phi}=\sum p^j \phi^j$ as a probability measure that assigns certain masses $\{p^j\}$ on an at most countable subset $\{\phi^j\}\subset\Phi$. Denote by $\bar{\Phi}$ the set of all quantizers, deterministic or random. Note that a deterministic quantizer can be thought of as a randomized one that assigns probability one to itself. For a randomized quantizer $\bar{\phi}=\sum p^j \phi^j$,  define its K-L divergences as the weighted average of those of the deterministic ones it randomizes:
\[I(m,m';{\bar{\phi}})=\sum p^j I(m,m';{\phi^j}).\]
This divergence for randomized quantizer will be key to our theorems.

The following assumption ensures basic regularities of the pdf's, it will be imposed throughout the rest of the paper.

\begin{assumption}\label{ass:1stmmts}
  For any two states $0\le m\neq m'\le M-1$, 
  \begin{displaymath}
	\textbf{E}_m\left[ \log\frac{f_m(X)}{f_{m'}(X)} \right]<\infty.
  \end{displaymath}
\end{assumption}

\section{Our Proposed Test $\delta_{A}(c)$}\label{sec:ABD}

In this section we will use tandem quantizers to define a class of ``two-stage'' tests $\delta(c)$, and show that asymptotic Bayes tests  can be found within it. The intuition is that the fusion center first makes a guess about the true state of nature and then tries to optimize the test based on the guess. 

As discussed in Mei \cite{mei}, tandem quantizer denotes the case when each sensor has the choice between two different sensor quantizers with at most one switch between them. Obviously, a tandem quantizer is the simplest non-stationary quantizers from the viewpoint of the number of switches. For the purpose of defining the two-stage sequential test $\delta(c)$, a useful alternative way to think about tandem quantizers is to divide the decision making into two stages.

In the \textit{first stage} of $\delta(c),$ one can use whatever reasonable stationary quantizers to make a preliminary decision
on which of the $M$ hypotheses is likely true, and the only requirement is that the sample size of this stage is large but is small relative to the overall sample sizes (or that of the second stage). Specifically, as $c \rightarrow 0,$ consider a sequence of $u(c) \in (0,1/2)$ such that $u(c)\to 0$ and $\log u(c)/\log c\to 0$,
e.g., $u(c)=1/|\log c|,$ and assume there is a  quantizer $\phi^0\in\Phi$ such that for any $0\le m,m'\le M-1$,
\begin{equation}\label{equ:prmlqntzrDecnt}
  I(m,m';{\phi^0})>0.
\end{equation}
Now in the first stage, the local sensor uses the stationary quantizer $\phi^0$ to send sensor messages to the fusion center, which will then face the classical multihypothesis sequential detection problem based on the i.i.d. quantized sensor messages $\phi^0(X_n).$ Hence, one can  recursively update the posterior distribution $(\pi_{0,n},\dots,\pi_{M-1,n}),\quad n=1,2,\dots$  at the fusion center as follows:
\begin{equation} \label{equ:postupdate}
  \pi_{m,n}=\frac{\pi_{m,n-1}f_m(U_n; \phi^0)}{\sum_{0\le m'\le M-1}\pi_{m',n-1}f_{m'}(U_n; \phi^0)},
\end{equation}
where $U_n$ is the quantized message at time $n$. As a reasonable test for the preliminary decision, the fusion center will stop the first stage at time $N^0$:
\begin{equation*}\label{equ:EndOfStage1}
  N^0=\min\{n\ge 1: \max_{0\le m\le M-1}{\{\pi_{m,n}\}}\ge 1-u(c)\},
\end{equation*}
and decides that the preliminary decision $D^0$ of the most promising state of nature is
\begin{equation*}\label{equ:dfnd0predcsn}
D^0=\argmax_{0\le m\le M-1}\pi_{m,N^0}.
\end{equation*}

In the \textit{second stage} of our two-stage test $\delta(c),$ the local sensor switches to another stationary (though likely randomized) quantizer, whose choices will likely depend on the  preliminary decision $D^{0}$ of the first stage.  Denote the quantizer used in the second stage as $\bar{\phi}_m$ when $D^0=m$, where $m=0,1,\dots, M-1$.

In the second stage, with the new quantizer applied at the local sensor, the fusion center starts afresh to update the posterior distribution $(\pi_{0,n},\dots,\pi_{M-1,n})$ based on i.i.d. sensor messages in the second stage. An efficient stopping rule for the fusion center can then be found as in Dragalin et al. \cite{dtv} as follows. Let $ r_{m,n}=\sum_{m'\neq m} \pi_{m',n} W(m',m)$ be the average loss by making a decision $m$ at time $n$, and let $r'_{m,n}=\min_{m'\neq m} \pi_{m,n}W(m,m')$ be the least value of loss by making some decision $m'\neq m$ at time $n$ while $m$ is the true state of nature. Define a total of $M$ stopping times:
\begin{equation}\label{equ:timeToStop2Stage}
N_m=\{n\ge N^0: \frac{r'_{m,n}}{r_{m,n}}>\frac{1}{c}\}, \quad m=0,1,\cdots, M-1.
\end{equation}
The fusion center can stop the second stage (hence the whole procedure) at time $N=\min\{N_m: 0\le m\le M-1\}$, and makes a final decision $D=m$ if $N=N_m$.

It is worth discussing the implementation of the likely randomized quantizer $\bar{\phi}_{m}$ if $D^0=m$ is the preliminary decision. We also need to give a explicit formula for updating posterior when randomized quantizer is used to form reports. Suppose $\bar{\phi}_{m}=\sum p^j \phi^j$. The key of any allowable randomization schemes is that the fusion center must know which deterministic quantizer is finally chosen, otherwise it may lose significant information and compromise the decision making efficiency.  We propose two alternative ways to achieve this goal. The most straightforward way is to let the fusion center do the randomization directly. Specifically, at a time step $n$ of the second stage, the fusion center selects a deterministic quantizer $\phi^{j}$ randomly according to the probability measure $\{p^j\},$ and informs the local sensor its choice through a feedback. Meanwhile, the posteior distributions should be updated as follows:
\begin{equation}\label{equ:postupdaterdmnqntzr}
  \pi_{m,n}=\frac{\pi_{m,n-1}f_m(U_n; \phi^j)}{\sum_{0\le m'\le M-1}\pi_{m',n-1}f_{m'}(U_n; \phi^j)}.
\end{equation}
An alternative way of randomization is to implement a ``block design'' at local level. Suppose that $\bar{\phi}_m$ is randomized by a finite number, say $J,$ of deterministic quantizers, and $b$ is a common denominators of the rational probabilities $p^1, \ldots, p^{J}$. Then take ``blocks" of $b$ observations, and in each block $\phi^1$,\dots, $\phi^J$ are used following a fixed order such that each $\phi^j$ appears exactly $b p^j$ times. In this way the fusion center also knows which quantizer is used at each time step and it will update the posterior just as in (\ref{equ:postupdaterdmnqntzr}).

For our proposed two-stage procedure $\delta(c),$ its asymptotic properties are summarized in the following theorem, whose proof is omitted since it can be derived along the same lines as those in Section V of Kiefer and Sacks \cite{ks}. To state the theorem, first we define the following information number for a quantizer $\bar{\phi}\in\bar{\Phi}$ and state $m=0,1,\dots, M-1$:
\begin{equation} \label{Jan11eq1}
  I(m; \bar{\phi})=\min_{m'\neq m}I(m, m'; \bar{\phi}).
\end{equation}
\begin{theorem} \label{the:asymp2stageproc}
  Let $\{\bar{\phi}_m: m=0,1,\dots, M-1\}$ be the randomized quantizers applied in the second stage of $\delta(c)$, and each $\bar{\phi}_m$ randomizes finite number of deterministic quantizers. Suppose $I(m; \bar{\phi}_m)>0$, $\pi_m>0$ for any $m$. Then as $c\to 0$, for the sample size $N$:
  \begin{equation} \label{equ:samplesize2stageproc}
	\textbf{E}_m[N]=(1+o(1))|\log c|/I(m; \bar{\phi}_m), \quad m=0,1,\dots,M-1,
  \end{equation}
  and for the probability of incorrect decisions:
  \begin{equation}
	\label{equ:wrongdecision2stageproc}
	\textbf{P}_m[D\neq m]=O(c), \quad m=0,1,\dots, M-1.
  \end{equation}
  Thus, the Bayes risk of the proposed two-stage test $\delta(c)$ is given by
  \begin{equation}
	\label{equ:risk2stageproc}
	\mathcal{R}_c(\delta)=c|\log c|(1+o(1)) \sum_{0\le m\le M-1} \frac{\pi_m}{I(m; \bar{\phi}_m)}.
  \end{equation}
\end{theorem}

In light of Theorem \ref{the:asymp2stageproc}, from the asymptotic viewpoint, an optimal procedure within the class of two-stage tests should maximize the information numbers $I(m; \bar{\phi}_m)$ so as to minimize the Bayes risk. This leads to a natural definition of the optimal quantizers that we should use in the second stage:

\begin{definition}
For $m=0,1,\ldots, M-1,$ the quantizer $\bar{\phi}^{\textrm{max}}_m$ is defined as the maximin quantizers with respect to $\textbf{P}_m$ if
\begin{equation*}\label{equ:maximins}
  \bar{\phi}^{\textrm{max}}_{m}=\argsup_{\bar{\phi}\in\bar{\Phi}}(I(m;{\bar{\phi}})).
\end{equation*}
\end{definition}
                                                                                Let us focus on the two-stage procedure $\delta_A(c)$ with the maximin quantizers being applied on the second stage. In next section, we will show that each $\bar{\phi}^{\textrm{max}}_m$ can be attained by randomizing at most $M-1$ deterministic quantizers. Hence by Theorem \ref{the:asymp2stageproc}, it has a Bayes risk
 \begin{equation}\label{equ:BayesOptimal}
\mathcal{R}_c(\delta_A(c))=(1+o(1)) c|\log c| \sum_m \frac{\pi_m}{I(m)}
  \end{equation}
as $c\to 0$, where $I(m)=\sup_{\bar{\phi}\in\bar{\Phi}} I(m;{\bar{\phi}})$.

Surprisingly, test $\delta_A(c)$ is not only the best among the two-stage tests, but also an asymptotically Bayes solution to problem (P2). This is a direct consequence of the following important theorem:

\begin{theorem} \label{the:asympoptdeltaIc}
Relation (\ref{equ:BayesOptimal}) is also satisfied by $\delta_B^*(c)$, the Bayes procedure.  
\end{theorem}
\proof The conclusion will be established once we prove the following: for any test with the probability of making incorrect decisions $\textbf{P}_m(D\neq m)=O(c\log c)$ for $m=0,1,\dots, M-1$, its expected values of the total time steps must satisfy $\textbf{E}_m N\ge (1+o(1))|\log c|I(m)^{-1}$ for any state $m$ as $c\to 0$. However this can be proved in the same way as Theorem 1 of Tsitovich \cite{tsito}. 
\endproof
\medskip
It is useful to point out that although the stopping rules of the asymptotic Bayes test $\delta_A(c)$  involve the prior distribution $\{\pi_m\}$'s, this is not essential and the key is for the local sensor to use the maximin quantizers $\bar{\phi}^{\textrm{max}}_m$'s at the second stage. In fact, since the maximin quantizers does not depend on the prior distribution $\{\pi_m\}$'s, (\ref{equ:samplesize2stageproc}) and (\ref{equ:wrongdecision2stageproc}) show that the optimality of $\delta_A(c)$ is robust w.r.t. a priori distribution $\{\pi_m\}$ as long as its support covers all $M$ possible states of nature. 

\section{Characterizing the Maximin Quantizers}\label{sec:SearchMaximinQntzr}

In this section, we provide a deeper understanding of the maximin quantizers $\{\bar{\phi}^{\textrm{max}}_m: m=0,1,\dots,M-1\}$ and also illustrate how to compute them explicitly when the sensor messages are binary. For this purpose, we first introduce the concept of
the unambiguous likelihood quantizer (ULQ), which was proposed in Tsitsiklis \cite{tsi} as a generalization of Monotone Likelihood Ratio Quantizer (MLRQ).

For simplicity, we assume that for any set of real numbers $\{a_{m'}: 0\le m'\le M-1\}$ which are not all zeros, 
\begin{equation}\label{equ:ULQIsExtCondition}
  P_{m}(\sum_{m'} a_{m'} f_{m'}(X)=0)=0, \quad 0\le m\le M-1.
\end{equation}
Note that (\ref{equ:ULQIsExtCondition}) is easily satisfied by the common continuous pdf families like normal, exponential, etc. 

\begin{definition}\label{def:ULQ}
  Under (\ref{equ:ULQIsExtCondition}), a deterministic quantizer $\phi\in\Phi$ is said to be an unambiguous likelihood quantizer if there exist real numbers $\{a_m: 0\le m\le M-1\}$ which are not all zero, such that 
\begin{equation*}\label{equ:dfnULQ1}
\phi(X)=I(\sum_m a_m f_m(X)>0),
\end{equation*}
\end{definition}
It is easy to see that in the case of binary simple hypothesis testing, i.e., $M=2$, the ULQs become MLRQs.

With the definition of ULQs, now it is time to state the following useful theorem which characterizes the maximin quantizers $\{\bar{\phi}^{\textrm{max}}_m\}$.
\begin{theorem}\label{the:MaximinQntzrApproxByULQ}
  Under (\ref{equ:ULQIsExtCondition}), each maximin quantizer $\bar{\phi}^{\textrm{max}}_m$ can be attained as a randomization of at most $M-1$ ULQs.
\end{theorem}

The detailed proof involves tedious technical details, and thus here we will only provide a high-level short explanation. For a fixed state $m$, finding the maximin quantizers against the other $M-1$ states is equivalent to solving an optimization problem in an $M-1$ dimensional space, where each quantizer, deterministic or randomized, corresponds to a point in it. By Tsitsiklis \cite{tsi}, these points construct a convex region whose extremal points all correspond to ULQs under the condition of (\ref{equ:ULQIsExtCondition}). Moreover,  the maximin quantizers correspond to the points that must be on the surface of the convex region, and thus can be expressed as a convex combination of at most $M-1$ extremal points (see Hormander \cite{hor}). Combining these results together leads to the desired relation between the  maximin quantizers and the ULQs.
\medskip

With Theorem \ref{the:MaximinQntzrApproxByULQ}, we are ready to illustrate how to find the maximin quantizers numerically. 

Fix any state $m$, define $M^2-1$ parameters as probability masses
$\{p^j_m: 1\le j\le M-1, p^j_m\ge 0, \sum_j p^j_m=1\},$
and ULQ coefficients
$\{a_{m,m'}^{j}: 1\le j\le M-1, 0\le m'\le M-1, \sum_{m'}  (a_{m,m'}^{j})^2=1\}.$
Based on every combination of these parameters, define by $\bar{\phi}$ the quantizer randomizing $M-1$ ULQs:
$\bar{\phi}=\sum_{j=1}^{M-1} p_{m}^j \phi_{m}^j,$
where
\begin{equation*}\label{equ:ULQDcmpdFmMaximinForm}
  \phi_{m}^j(X)=I(\sum_{m'} a^j_{m,m'} f_{m'}(X)>0).
\end{equation*}
The maximin quantizer $\bar{\phi}^{\textrm{max}}_m$ can then be found as $\bar{\phi}$ that maximizes
\begin{equation}\label{equ:ObjFunc}
\min_{l \neq m} I(m,l;{\bar{\phi}}),
\end{equation}
among all possible combinations of $\{p^j_m; a^j_{m,m'}\}.$

\section{Examples}
In this section we illustrate our procedure with a concrete example. Suppose that the raw sensor observations $X_n$'s are distributed according to $N(\mu,1).$ If there are only $M=2$ hypotheses on $\mu,$ say testing $\textbf{H}_0: \mu = 0$ against $\textbf{H}_1: \mu = 1,$ then
there is no randomization involved in the second stage, and the maximin quantizer is just the ULQs which becomes the MLRQs when $M=2.$ Such a result is consistent with those in Mei \cite{mei}.

Now suppose there are $M=3$ hypotheses regarding the normal mean:
$\textbf{H}_0: \mu=0,$ $\textbf{H}_1: \mu=-1,$ and $\textbf{H}_2: \mu=1.$
For this specific case, it is not too difficult to solve the optimization problem (\ref{equ:ObjFunc}) by linear programming. Up to the precision of four decimal places, numerical computations show that all three maximin quantizers turn out to be deterministic ones:
$\phi_0 =  I(X>0),$ $\phi_1 = I(X>-0.7941),$ and $\phi_2 =  I(X>0.7941),$
and their corresponding maximin information numbers are
$I(0)=0.3137$ and $I(1)=I(2)=0.3186$. 
For the first stage, the quantizer $\phi_0$ can be applied because it satisfies the condition (\ref{equ:prmlqntzrDecnt}). By Theorem \ref{the:asymp2stageproc}, the risk of $\delta_A(c)$ can be approximated by
\begin{displaymath}
  \mathcal{R}_c(\delta_A(c))=c|\log c|(1+o(1))(\frac{\pi_0}{0.3137}+\frac{\pi_1+\pi_2}{0.3186}).
\end{displaymath}

As a comparison, in the centralized version when the whole raw observations are allowed to be used at the fusion center, it can be shown that the Bayes risk of the optimal centralized test is
\begin{displaymath}
  \mathcal{R}_c(\delta_{\textrm{cen}}(c))=2 c|\log c| (1+o(1)),
\end{displaymath}
see, for example,  Dragalin et al. \cite{dtv} and Kiefer and Sacks \cite{ks}. Thus the asymptotic efficiency of $\delta_A(c)$ with respect to the optimal centralized test is
\[
\lim_{c \rightarrow 0} \mathcal{R}_c(\delta_{\textrm{cen}}(c)) / \mathcal{R}_c(\delta_A(c)) \ge 2 / (1/0.3137) = 0.6274.
\]
In particular, if we just merely introduce another identical sensor into the network system, then the efficiency of $\delta_A(c)$ will be doubled and the corresponding decentralized test will have better properties than that of $\delta_{\textrm{cen}}(c)$.

\section{Conclusion}

In this article, the problem of decentralized testing multihypotheses in (single) sensor networks is studied. Asymptotically Bayes test $\{\delta_A(c)\}$  is constructed by combining the ideas of ``tandem quantizers'', ``unambiguous quantizers'', and ``randomized quantizers.'' Such a test involves a new concept of maximin quantizers which are discussed in details, both theoretically and numerically.

It is natural to extend our results to the networks with multiple sensors, where different sensors may use different quantizers. A more interesting extension is to understand what happens when one or more hypotheses are not simple, i.e., the composite multihypotheses case. These will be reported elsewhere.



\section*{Acknowledgment}
This work was supported in part by the AFOSR grant FA9550-08-1-0376 and the NSF Grant CCF-0830472.

\end{document}